\documentclass[11pt]{article}
\usepackage{amssymb,latexsym}
\usepackage{epsfig}
\usepackage{eufrak}
\usepackage{amsmath}
\usepackage{mathrsfs}
\usepackage{color}

\newtheorem{theorem}{Theorem}
\newtheorem{lemma}{Lemma}

\newcommand{\be}{\begin{equation}}
\newcommand{\ee}{\end{equation}}
\newcommand{\bee}{\begin{eqnarray*}}
\newcommand{\eee}{\end{eqnarray*}}
\newcommand{\bel}{\begin{eqnarray}}
\newcommand{\eel}{\end{eqnarray}}
\newcommand{\bec}{\begin{cases}}
\newcommand{\eec}{\end{cases}}
\newcommand{\bem}{\begin{bmatrix}}
\newcommand{\eem}{\end{bmatrix}}

\newcommand{\la}{\label}
\newcommand{\li}{\left}
\newcommand{\ri}{\right}

\newcommand{\ovl}{\overline}
\newcommand{\udl}{\underline}

\newcommand{\lc}{\lceil}
\newcommand{\rc}{\rceil}
\newcommand{\lf}{\lfloor}
\newcommand{\rf}{\rfloor}

\newcommand{\vep}{\varepsilon}

\newcommand{\de}{\delta}

\newcommand{\vDe}{\varDelta}

\newcommand{\se}{\theta}

\newcommand{\ze}{\zeta}
\newcommand{\al}{\alpha}
\newcommand{\ba}{\beta}

\newcommand{\Om}{\Omega}

\newcommand{\f}{\frac}
\newcommand{\sq}{\sqrt}
\newcommand{\cd}{\cdots}

\newcommand{\qu}{\quad}
\newcommand{\qqu}{\qquad}

\newcommand{\mscr}{\mathscr}

\newcommand{\bb}{\mathbb}

\newcommand{\wh}{\widehat}

\newcommand{\mrm}{\mathrm}
\newcommand{\bs}{\boldsymbol}

\newcommand{\ap}{\approx}

\newcommand{\LRA}{\Longleftrightarrow}
\newcommand{\sh}{\slash}

\newcommand{\tx}{\text}

\newcommand{\iy}{\infty}

\newcommand{\pa}{\partial}

\newcommand{\bed}{\begin{description}}
\newcommand{\eed}{\end{description}}
\newcommand{\bei}{\begin{itemize}}
\newcommand{\eei}{\end{itemize}}
\newcommand{\ben}{\begin{enumerate}}
\newcommand{\een}{\end{enumerate}}

\newcommand{\beL}{\begin{lemma}}
\newcommand{\eeL}{\end{lemma}}
\newcommand{\beT}{\begin{theorem}}
\newcommand{\eeT}{\end{theorem}}
\newcommand{\sect}{\section}

\newcommand{\bpf}{\begin{pf}}
\newcommand{\epf}{\end{pf}}
\newcommand{\bsk}{\bigskip}
\newcommand{\bi}{\binom}

\newcommand{\pfbox}{\hfill\mbox{$\Box$}}

\newenvironment{pf}{\paragraph*{Proof{\rm.}}}{\pfbox\bigskip}

\begin{document}

\title{{\bf Exact Computation of Minimum Sample Size for Estimation of  Binomial Parameters}
\thanks{The author had been previously working with Louisiana State University at Baton Rouge, LA 70803, USA,
and is now with Department of Electrical Engineering, Southern University and A\&M College, Baton Rouge, LA
70813, USA; Email: chenxinjia@gmail.com}}

\author{Xinjia Chen}

\date{July 2007}

\maketitle

\begin{abstract}

It is a common contention that it is an ``impossible mission'' to exactly determine the minimum sample size for
the estimation of a binomial parameter with prescribed margin of error and confidence level. In this paper, we
investigate such a very old but also extremely important problem and demonstrate that the difficulty for
obtaining the exact solution is not insurmountable.  Unlike the classical approximate sample size method based
on the central limit theorem, we develop a new approach for computing the minimum sample size that does not
require any approximation. Moreover, our approach overcomes the conservatism of existing rigorous sample size
methods derived from Bernoulli's theorem or Chernoff bounds.

Our computational machinery consists of two essential ingredients. First, we prove that the minimum of coverage
probability with respect to a binomial parameter bounded in an interval is attained at a discrete set of finite
many values of the binomial parameter.  This allows for reducing infinite many evaluations of coverage
probability to finite many evaluations. Second, a recursive bounding technique is developed to further improve
the efficiency of computation.

\end{abstract}

\sect{Introduction}

The estimation of a binomial parameter is a fundamental problem in probability and statistics.  The practical
importance of such estimation problem can be seen by its numerous applications in various fields of sciences and
engineering.  Specifically, the problem is formulated as follows.

Let $X$ be a Bernoulli random variable defined in a probability space $(\Om, \mscr{F}, \Pr)$ such that $\Pr \{ X
= 1 \} = p$ and $\Pr \{ X = 0 \} = 1 - p$ with $p \in (0,1)$.  It is a frequent problem to estimate $p$ based on
$n$ identical and independent samples $X_1, \cd, X_n$ of $X$.  The parameter $p$ is referred to as a binomial
parameter, since it defines a binomial experiment for a given sample size $n$.

An estimate of $p$ is conventionally taken as $\wh{\bs{p}}_n = \f{ \sum_{i=1}^n X_i } { n }$.  The nice property
of such estimate is that it is of maximum likely-hood and possesses minimum variance among all unbiased
estimates.  A crucial question in the estimation is as follows:

\bsk

{\it Given the knowledge that $p$ belongs to interval $[a, b]$, what is the minimum sample size $n$ that
guarantees the difference between $\wh{\bs{p}}_n$ and $p$ be bounded within some prescribed margin of error with
a confidence level higher than a prescribed value?}

\bsk

It is generally believed that an exact answer to this fundamental question is not possible (see, e.g.,
\cite{Rahme} and the references therein).  However, our recent investigation shows that an exact solution can be
found by combining the power of mathematical analysis and modern computers.  The main contribution of this paper
is to provide an exact answer to this important question. In contrast to existing methods in the literature, we
aim at finding rigorous solutions while avoiding unnecessary conservatism.

The paper is organized as follows. In Section 2, the techniques for computing the minimum sample size is
developed with the margin of error taken as a bound of absolute error. In Section 3, we derive corresponding
sample size method by using relative error bound as the margin of error.  In Section 4, we develop techniques
for computing minimum sample size with a mixed error criterion.  Section 5 is the conclusion. The proofs are
given in Appendices.

Throughout this paper, we shall use the following notations. The set of integers is denoted by $\bb{Z}$.  The
ceiling function and floor function are denoted respectively by $\lc . \rc$ and $\lf . \rf$ (i.e., $\lc x \rc$
represents the smallest integer no less than $x$; $\lf x \rf$ represents the largest integer no greater than
$x$).  For non-negative integer $m$, the combinatoric function $\bi{m}{z}$ with respect to integer $z$ means
\[
\bi{m}{z} = \bec \f{ m! } { z! (m- z)! } & \tx{for} \; 0 \leq z \leq m,\\
0 & \tx{for} \; z < 0 \; \tx{or} \; z > m. \eec
\]
The binomial function is $B(n, k, p) = \bi{n}{k} p^k (1 - p)^{n-k}$, assuming value $0$ for $k < 0$ or $k > n$.
The summation of binomial function is denoted as {\small $S(n, k, l, p) = \sum_{i=k}^l B(n, i, p)$}.  The left
limit as $\eta$ tends to $0$ is denoted as $\lim_{\eta \downarrow 0}$. The notation ``$\LRA$'' means ``if and
only if''. The other notations will be made clear as we proceed.

\section{Control of Absolute Error}
In this section, we shall review the classical sample size problem and the existing solutions. In particular, we
shall elaborate the difficulty that has been considered as insurmountable in the literature. We will demonstrate
that such ``seemingly'' insurmountable difficulty can be made disappear by a careful analysis of the coverage
probability.

Formally, the classical sample size problem is stated as follows.  Let $\vep \in (0,1)$ be the margin of
absolute error and $\de \in (0,1)$ be the confidence parameter. In many applications, it is desirable to find
the smallest sample size $n$ such that \be \la{absn}
 \Pr \li \{ |\wh{\bs{p}}_n - p|
< \vep \ri \} > 1 - \de \ee for any $p \in [a,b]$.  Here $\Pr \li \{ |\wh{\bs{p}}_n - p| < \vep \ri \}$ is
referred to as the coverage probability.  The interval $[a, b]$ is introduced to take into account the knowledge
of $p$. If no knowledge of $p$ is available, the interval $[a, b]$ can be taken as $[0,1]$.

The classical sample size problem associated with (\ref{absn}) has been extensively studied in the literature.
As pointed out in pages 83-84 of \cite{Rahme}, it is commonly believed that the exact computation of the minimum
sample size is impossible.   This is due to the intuitive that, for suitably chosen $k_1$ and $k_2$,  \be
\la{dir}
 \Pr \li \{ |\wh{\bs{p}}_n - p| < \vep \ri \} = \sum_{k = k_1}^{k_2} \bi{n}{k} p^k (1 - p)^{n - k},
\ee where both the summand and $k_1, \; k_2$ depend on the unknown value $p$, making the direct use of
(\ref{dir}) and therefore (\ref{absn}) ``almost impossible in practice.'' Such argument is very typical and can
be seen in page 84, lines 1-6 of \cite{Rahme}.   In general,  one tends to think that infinite many evaluations
of the right-hand side of (\ref{dir}) is required to determine whether the coverage probability is greater than
$1 - \de$ for any $p$ in interval $[a, b]$. Motivated by the ``seemingly'' prohibitive computational complexity,
statisticians have been settled to finding approximation or conservative bounds for the minimum sample size
associated with (\ref{absn}).

The conventional solution is based on the normal approximation (see, e.g., \cite{Desu, Fishman} and the
references therein).  The drawback of such sample size method is that the coverage probability in (\ref{absn})
may be significantly below the prescribed confidence level $1 - \de$.  This can be an extremely severe problem
in the case that the upper bound, $b$,  of the binomial parameter is small.    Such criticism is very usual as
can be seen in \cite{Fishman, Hampel} and many other literatures.  The issue of inaccuracy remains significant
even for the case that no information of $p$ is available, i.e., $[a, b] = [0,1]$.  In this case, the minimum
sample size is approximated as \be \la{norm} n_{\mrm{G}} \ap  \f{ Z_{\de \sh 2}^2 } { 4 \vep^2 }  \ee where
$Z_{\de \sh 2}$ satisfies $\int_{ Z_{\de \sh 2} }^\iy \f{1}{\sq{2 \pi}} e^{ - \f{x^2}{2}  } d x = \f{\de}{2}$.
Application of the approximate formula  (\ref{norm}) must introduce unknown error in reporting the statistical
accuracy of the estimation of $p$. In order to eliminate the inaccuracy of normal approximation, one can resort
to the large deviation type inequalities to derive an upper bound for the minimum sample size.  A well-known
result is the Chernoff bound, which asserts that (\ref{absn}) is true  for any $p \in [0,1]$ provided that
 \be
 \la{Chern}
 n > \f{ \ln \f{2}{\de} } { 2 \vep^2 }.
 \ee
The Chernoff bound significantly improves upon the sample size bound provided by the famous Bernoulli theorem,
which states that (\ref{absn}) is ensured  for any $p \in [0,1]$ if \be \la{ber} n > \f{1}{4 \vep^2 \de}. \ee
The major problem of sample size formulas (\ref{Chern}) and (\ref{ber}) is the unduly conservativeness.  The
sample size obtained from (\ref{Chern}) or (\ref{ber}) can be substantially larger than the minimum sample size.

\bsk

Since one of the fundamental goals of statistics is to provide rigorous and the least conservative
quantification of uncertainty in statistical inference, it is a persistent concern of statisticians and
practitioners to determine the exact value of minimum sample size associated with (\ref{absn}).  After a
thorough investigation, we discovered that the exact determination of minimum sample size is readily tractable
with modern computational power by taking advantage of the behavior of the coverage probability characterized by
Theorem \ref{thm_abs} as follows.

\beT \la{thm_abs} Let $0 < \vep < 1$ and $0 \leq a < b \leq 1$.  Let $X_1, \cd, X_n$ be identical and
independent Bernoulli random variables such that, for $i = 1, \cd, n$,  $\Pr \{ X_i = 1 \} = 1 - \Pr \{ X_i = 0
\} = p$ with $p \in [a, b]$. Let $ \wh{\bs{p}}_n = \f{ \sum_{i = 1}^n X_i  }{ n }$.
Then, the minimum of $\Pr \{
| \wh{\bs{p}}_n - p | < \vep \}$ with respect to $p \in [a, b]$ is achieved at the finite set $\{a, b \} \cup
 \{ \f{\ell}{n} + \vep \in (a, b) : \ell \in \bb{Z} \} \cup
 \{ \f{\ell}{n} - \vep \in (a, b) : \ell \in \bb{Z} \} $,
 which has less than $2n(b-a) + 4$ elements.
\eeT

See Appendix A for a proof.  The application of Theorem \ref{thm_abs} in the computation of minimum sample size
is obvious.  For a fixed sample size $n$, since the minimum of coverage probability with $p \in [a, b]$ is
attained at a finite set, it can determined by a computer whether the sample size $n$ is large enough to ensure
(\ref{absn}) for any $p \in [a, b]$. Starting from $n = 2$,  one can find the minimum sample size by gradually
incrementing $n$ and checking whether $n$ is large enough.

By the fact of symmetry that $\Pr \{ | (1 - \wh{\bs{p}}_n) - (1 - p)
| < \vep \} = \Pr \{ | \wh{\bs{p}}_n - p | < \vep \}$, we can
restrict $p$ to a smaller range $[a^\prime, b^\prime]$ such that
\[
\min_{p \in [a, b ]} \Pr \{ | \wh{\bs{p}}_n - p | < \vep \} = \min_{p \in [a^\prime, b^\prime ]} \Pr \{ |
\wh{\bs{p}}_n - p | < \vep \}
\]
where
\[
a^\prime = \bec
 a & \tx{for} \; a + b \leq 1,\\
 1 - b & \tx{for} \; a + b > 1
\eec \qqu  b^\prime = \bec
 b & \tx{for} \; b \leq \f{1}{2},\\
 \f{1}{2} & \tx{for} \; a < \f{1}{2} < b,\\
 1 - a & \tx{for} \; a \geq \f{1}{2}.
\eec
\]
Clearly, $0 \leq a^\prime < b^\prime \leq \f{1}{2}$ and $b^\prime - a^\prime \leq b - a$.  Hence, without loss
of any generality, we can assume $0 \leq a < p < b \leq \f{1}{2}$.

As can be seen from Theorem \ref{thm_abs}, for $a = 0$ and $b =
\f{1}{2}$, the total number of binomial summations to be evaluated
is no more than $n + 2$, since the coverage probability for $a =
0$ is one and no computation is needed.

For computational purpose, we have

\beT \la{thm_com} Let $0 \leq a < b \leq \f{1}{2}$ and $0 < \vep <
\f{1}{2}$. Define {\small \bee \mscr{S} & = &   \li \{ c_{+}
(\ell) \mid 1 + \lf n(a - \vep) \rf \leq \ell \leq \lc n(b - \vep)
\rc - 1 \ri \} \cup \li \{c_a, \; c_b
\ri \}\\
&  & \cup \li \{ c_{-} (\ell) \mid 1 + \lf n(a + \vep) \rf  \leq
\ell \leq \lc n(b + \vep) \rc -  1 \ri \}  \eee} where

$c_a =  \sum_{k = \lf n(a - \vep) \rf + 1 }^{ \lc n(a + \vep) \rc
-1 } B(n, k, a), \qqu \qu c_b = \sum_{k = \lf n(b - \vep) \rf + 1
}^{ \lc n(b + \vep) \rc -1 } B(n, k, b)$,

$c_{+} (\ell)  =  \sum_{ k = \ell + 1 }^{ \ell - 1 + \li \lc 2 n
\vep \ri \rc  } B \li (n, k,  \f{ \ell } { n } + \vep \ri ), \qqu
c_{-}(\ell)  =  \sum_{ k = \ell + 1 - \li \lc 2 n \vep \ri \rc }^{
\ell  - 1 } B \li (n, k, \f{ \ell } { n } - \vep \ri )$ for $\ell
\in \bb{Z}$.

Define {\small \bee &
  & \ovl{\vDe} (n,\se, r,s)
=   B(n-1, r-1, \se) +  B \li (n, s + 1, \se + \f{1}{n} \ri ) -  B(n, r-1, \se)  - B(n-1, s, \se)\\
&   &  \qu +  \; \f{ n-1 } {2 n } \li [ \ovl{B}(n-2,r-2, \se) +
\ovl{B}(n-2, s, \se) -   \udl{B}(n-2,r-1, \se)  - \udl{B}(n-2,
s-1, \se) \ri ],\\
&  & \udl{\vDe} (n,\se, r,s)
=   B(n-1, r-1, \se) +  B \li (n, s + 1, \se + \f{1}{n} \ri ) -  B(n, r-1, \se)  - B(n-1, s, \se)\\
&   &  \qu + \; \f{ n-1 } {2 n } \li [ \udl{B}(n-2,r-2, \se) +
\udl{B}(n-2, s, \se) -   \ovl{B}(n-2,r-1, \se)  - \ovl{B}(n-2,
s-1, \se) \ri ] \eee} with $\udl{B} (n, k, \se) = \min \li \{ B(n,
k, \se), \;  B \li ( n, k, \se+ \f{1}{n} \ri ) \ri \}$ and
\[ \ovl{B} (n, k, \se) =  \bec \max \li \{ B(n, k, \se), \; B
\li ( n, k, \se+ \f{1}{n} \ri
) \ri \}  & \tx{for} \; \f{k}{n} \notin [\se, \se + \f{1}{n}],\\
B \li ( n, k, \f{k}{n} \ri )  & \tx{for} \; \f{k}{n} \in [\se, \se
+ \f{1}{n}]. \eec
\]
Then, the following statements hold true:

(I) The minimum of $\Pr \{ | \wh{\bs{p}}_n - p | < \vep \}$ with
respect to $p \in [a, b]$ equals the minimum of $\mscr{S}$, i.e.,
$\min_{p \in [a, b ]} \Pr \{ | \wh{\bs{p}}_n - p | < \vep \} =
\min \mscr{S}$.

(II) For $\lf n(a - \vep) \rf \leq \ell \leq \lc n(b - \vep) \rc - 1$, \be \la{recur}
 [1 - c_+(\ell)] +
\udl{\vDe} (n, \se_\ell, r_\ell, s_\ell) \leq 1 - c_+(\ell-1) \leq  [1 - c_+(\ell)] + \ovl{\vDe} (n, \se_\ell,
r_\ell, s_\ell), \ee where $\se_\ell = \f{ \ell - 1 } { n } + \vep ,  \;  r_\ell = \ell + 1$ and $s_\ell = \ell
- 2 + \li \lc 2 n \vep \ri \rc$.

(III) For $\lf n(a + \vep) \rf  \leq \ell \leq \lc n(b + \vep) \rc - 1$, \be \la{recur2}
 [1 - c_-(\ell)] +
\udl{\vDe} (n, \se_\ell^\prime, r_\ell^\prime, s_\ell^\prime) \leq 1 - c_-(\ell-1) \leq [1 - c_-(\ell)] +
\ovl{\vDe} (n, \se_\ell^\prime, r_\ell^\prime, s_\ell^\prime) \ee where $\se_\ell^\prime = \f{ \ell - 1 } { n }
- \vep , \; r_\ell^\prime = \ell + 1 - \li \lc 2 n \vep \ri \rc$ and $s_\ell^\prime = \ell - 2$.
 \eeT

See Appendix B for a proof.  \bsk

For the purpose of reducing roundoff error, we shall evaluate the complementary probability $1 - c_+(\ell), \; 1
- c_-(\ell), 1 - c_a, 1 - c_b$ and compare them with $\de$ to determine whether the sample size $n$ is large
enough.  Since the comparison between $1 - c_+(\ell)$ and $\de$ usually only requires bounds of $1 - c_+(\ell)$,
a large amount of computation can be saved if we start from $\ell = \lc n(b - \vep) \rc - 1$ and recursively
build the bounds of $1 - c_+(\ell - 1)$ from the bounds of $1 - c_+(\ell)$ by making use of (\ref{recur}).
Similarly, we can apply (\ref{recur2}) to reduce the computation for the comparison between $1 - c_-(\ell)$ and
$\de$. This computational technique is especially useful when the sample size is large.  Due to the recursive
nature, we call it as the {\it recursive bounding} technique.  It should be noted the bounds may become too
conservative to be useful as the number of recursive steps increases.  In that situation, the recursive process
should be restarted with exact computation for the current index $\ell$.

In order to reduce computational effort, the evaluation should be performed earlier for $\f{\ell}{n} \pm \vep$
closer to $\f{1}{2}$ for the purpose of earlier determination of whether the sample size is sufficiently large.
This computational trick is motivated by our computational experience that for many values of $\ell$, $1 -
c_+(\ell)$ is non-decreasing with respect to $\ell$.  The situation  is similar for $1 - c_-(\ell)$.

To demonstrate the feasibility of our computational method, we provide some sample size values in Table
\ref{table_abs} for the case that $[a, b] = [0,1]$, i.e., no information for $p$ is available. Actually, with a
few hours of computer running time, we have produced a MATLAB data file for a large number of combinations of
margin of absolute error and confidence level.  Although the computational complexity of our approach is much
higher than that of existing explicit formulas, the computer running time is not an issue since a large data
file of sample size can be created and saved for forever use.

\begin{table}
\caption{Table of Sample Sizes} \label{table_abs}
\begin{center}
\begin{tabular}{|c||c||c|||c||c||c|||c||c||c|}
\hline $\vep$ & $\de$ & $n$ & $\vep$ & $\de$ & $n$ & $\vep$ & $\de$
& $n$\\
\hline
\hline $0.1$ & $0.05$ & $101$ & $0.1$ & $0.01$ & $171$ & $0.1$ & $0.001$ & $276$\\
\hline $0.05$ & $0.05$ & $391$ & $0.05$ & $0.01$ & $671$ & $0.05$ & $0.001$ & $1091$\\
\hline $0.01$ & $0.05$ & $9651$ & $0.01$ & $0.01$ & $16601$ & $0.01$ & $0.001$ & $27101$\\
\hline $0.005$ & $0.05$ & $38501$ & $0.005$ & $0.01$ & $66401$ & $0.005$ & $0.001$ & $108301$\\
\hline $0.001$ & $0.05$ & $960501$ & $0.001$ & $0.01$ & $1659001$ & $0.001$ & $0.001$ & $2707001$\\
 \hline
\end{tabular}
\end{center}
\end{table}

\section{Control of Relative Error}

Let $\vep \in (0,1)$ be the margin of  relative error and $\de \in
(0,1)$ be the confidence parameter.  It is interesting to determine
the smallest sample size $n$ so that
\[
\Pr \li \{  \f{ |\wh{\bs{p}}_n - p| } {p} < \vep  \ri \}
> 1 - \de
\]
for any $p \in [a,b]$.  As has been pointed out in Section 2, an essential machinery is to reduce infinite many
evaluations of the coverage probability $\Pr \{ | \wh{\bs{p}}_n - p | < \vep p \}$ to finite many evaluations.
Such reduction can be accomplished by making use of Theorem \ref{thm_rev} as follows.

 \beT \la{thm_rev} Let $0 < \vep < 1$ and $0
< a  < b \leq 1$.  Let $X_1, \cd, X_n$ be identical and
independent Bernoulli random variables such that, for $i = 1, \cd,
n$,  $\Pr \{ X_i = 1 \} = 1 - \Pr \{ X_i = 0 \} = p$ with $p \in
[a, b]$. Let $ \wh{\bs{p}}_n = \f{ \sum_{i = 1}^n X_i  }{ n }$.
Then, the minimum of $\Pr \li \{  \f{ |\wh{\bs{p}}_n - p| } {p} <
\vep \ri \}$ with respect to $p \in [a, b]$ is achieved at the
finite set $\{a, b \} \cup
 \{ \f{\ell}{n(1 + \vep)} \in (a, b) : \ell \in \bb{Z} \} \cup
 \{ \f{\ell}{n (1 - \vep)} \in (a, b) : \ell \in \bb{Z} \}$,
 which has less than $2n(b-a) + 4$ elements.
\eeT

See Appendix C for a proof.  For computational convenience, we have

\beT \la{thm_rev_com} Let $0 \leq a < b \leq 1$ and $0 < \vep < 1$. Define {\small \bee \mscr{S}_{\mrm{rev}} & =
&   \li \{ c_{+} (\ell) \mid 1 + \lf n a(1 - \vep) \rf \leq \ell \leq \lc n b (1 - \vep) \rc - 1 \ri \} \cup \li
\{c_a, \; c_b
\ri \}\\
&  & \cup \li \{ c_{-} (\ell) \mid 1 + \lf n a (1 + \vep) \rf \leq
\ell \leq \lc n b(1 + \vep) \rc -  1 \ri \}  \eee} where

$c_a =  \sum_{k = \lf n a (1 - \vep) \rf + 1 }^{ \lc n a (1 +
\vep) \rc -1 } B(n, k, a), \qqu \qu c_b = \sum_{k = \lf n b (1 -
\vep) \rf + 1 }^{ \lc n b(1 + \vep) \rc -1 } B(n, k, b)$,

$c_{+} (\ell)  =  \sum_{ k = \ell + 1 }^{ \li \lc \f{ 1 + \vep}{1 - \vep } \ell \ri \rc  - 1} B \li (n, k,  \f{
\ell } { n (1 - \vep)} \ri ), \qqu c_{-}(\ell) = \sum_{ k = \li \lf \f{ 1 - \vep}{1 + \vep } \ell \ri \rf  + 1
}^{\ell - 1} B \li (n, k, \f{ \ell } { n (1 + \vep) } \ri )$ for $\ell \in \bb{Z}$.  Then, the minimum of $\Pr
\{ | \wh{\bs{p}}_n - p | < \vep p \}$ with respect to $p \in [a, b]$ equals the minimum of
$\mscr{S}_{\mrm{rev}}$, i.e., $\min_{p \in [a, b ]} \Pr \{ | \wh{\bs{p}}_n - p | < \vep p \} = \min
\mscr{S}_{\mrm{rev}}$. \eeT

Actually, a similar type of recursive bounding technique as Theorem \ref{thm_com} can be developed to improve
the computational efficiency.  Theorem \ref{thm_rev_com} can be proved by applying Theorem \ref{thm_rev} and the
following observations:

(i) $\Pr \{ | \wh{\bs{p}}_n - p | < \vep p \}$ assumes values $c_a$ and $c_b$ for $p = a$ and $b$ respectively.

(ii) For {\small $p = \f{\ell}{n (1 + \vep)} \in \{ \f{\ell}{n (1 + \vep) } \in (a, b) : \ell \in \bb{Z} \}$},
we have $\lc n p (1 + \vep) \rc - 1 = \ell - 1, \;\; \lf n p (1  - \vep) \rf + 1 = \li \lf \f{1 - \vep}{1 +
\vep} \ell \ri \rf + 1, \;\; c_-(\ell) = \Pr \{ | \wh{\bs{p}}_n - p | < \vep p \}$ and $1 + \lf n a (1 + \vep)
\rf \leq \ell \leq \lc n b (1 + \vep) \rc - 1$.

(iii) For {\small $p = \f{\ell}{n (1 - \vep)} \in \{ \f{\ell}{n (1 - \vep)} \in (a, b) : \ell \in \bb{Z} \}$},
we have $\lf n p (1 - \vep) \rf + 1 = \ell + 1$ and $\lc n p (1  + \vep) \rc - 1 = \li \lc \f{1 + \vep}{1 -
\vep} \ell \ri \rc - 1, \;\; c_+(\ell) = \Pr \{ | \wh{\bs{p}}_n - p | < \vep p \}$ and $1 + \lf n a (1 - \vep)
\rf \leq \ell \leq \lc n b (1 - \vep) \rc - 1$.

\section{Control of Absolute Error or Relative Error}

Let $\vep_a \in (0,1)$ and $\vep_r \in (0,1)$ be respectively the margins of absolute error and relative error.
Let $\de \in (0,1)$ be the confidence parameter.  In many situations, it is desirable to find the smallest
sample size $n$ such that \be \la{mixedd}
 \Pr \li \{ |\wh{\bs{p}}_n - p| < \vep_a \; \; \mrm{or} \;\;  \f{
|\wh{\bs{p}}_n - p| } {p} < \vep_r \ri \}
> 1 - \de \ee for any $p \in
[a,b]$.  To make it possible to compute exactly the minimum sample size associated with (\ref{mixedd}), we have
Theorem \ref{thm_abs_rev} as follows.

\beT \la{thm_abs_rev} Let $0 < \vep_a < 1, \; 0 < \vep_r < 1$ and $0 \leq a < \f{\vep_a}{\vep_r} < b \leq 1$.
 Let $X_1, \cd, X_n$ be identical and independent
Bernoulli random variables such that, for $i = 1, \cd, n$,  $\Pr \{ X_i = 1 \} = 1 - \Pr \{ X_i = 0 \} = p$ with
$p \in [a, b]$. Let $ \wh{\bs{p}}_n = \f{ \sum_{i = 1}^n X_i  }{ n }$. Then, the minimum of $\Pr \li \{
|\wh{\bs{p}}_n - p| < \vep_a \; \; \mrm{or} \;\; \f{ |\wh{\bs{p}}_n - p| } {p} < \vep_r  \ri \}$ with respect to
$p \in [a, b]$ is achieved at the finite set $\{a, b, \f{\vep_a}{\vep_r} \} \cup
 \{ \f{\ell}{n} + \vep_a \in (a, \f{\vep_a}{\vep_r}) : \ell \in \bb{Z} \} \cup
 \{ \f{\ell}{n} - \vep_a \in (\f{\vep_a}{\vep_r}, b) : \ell \in \bb{Z} \} \cup
 \{ \f{\ell}{n(1 + \vep_r)} \in (a, \f{\vep_a}{\vep_r}) : \ell \in \bb{Z} \} \cup
 \{ \f{\ell}{n (1 - \vep_r)} \in (\f{\vep_a}{\vep_r}, b) : \ell \in \bb{Z}
 \}$,
 which has less than $2n(b-a) + 7$ elements.
\eeT

\bsk

As can be seen from Theorem \ref{thm_abs_rev}, for $a = 0$ and $b = 1$, the total number of evaluations of
probability is no more than $2n + 4$.  The detailed proof of Theorem \ref{thm_abs_rev} is omitted since it can
deduced from Theorem \ref{thm_abs} and Theorem \ref{thm_rev} with the observation that
\[
\Pr \li \{ |\wh{\bs{p}}_n - p| < \vep_a \; \; \mrm{or} \;\;  \f{
|\wh{\bs{p}}_n - p| } {p} < \vep_r \ri \} = \bec \Pr \li \{
|\wh{\bs{p}}_n - p| < \vep_a  \ri \} & \tx{for} \; p \in \li [ a,
\f{\vep_a}{\vep_r} \ri ], \\
\Pr \li \{  \f{ |\wh{\bs{p}}_n - p| } {p} < \vep_r \ri \} & \tx{for}
\; p \in \li ( \f{\vep_a}{\vep_r}, b \ri ]. \eec
\]
Such observation also indicates that the sample size problem associated with (\ref{mixedd}) can be decomposed as
the sample size problems for the cases of absolute error and relative error discussed previously.

To show the effectiveness of our sample size method, we present some sample size numbers in Table
\ref{table_mixed} for the case that no information for the binomial parameter is available, i.e., $[a, b] =
[0,1]$. We would like to remark that the computation can be easily managed by any personal computer.

\begin{table}
\caption{Table of Sample Sizes ($\vep_r = 0.1$)} \label{table_mixed}
\begin{center}
\begin{tabular}{|c||c||c|||c||c||c|||c||c||c|}
\hline $\vep_a$ & $\de$ & $n$ & $\vep_a$ & $\de$ & $n$ & $\vep_a$ & $\de$
& $n$\\
\hline
\hline $0.05$ & $0.05$ & $391$ & $0.05$ & $0.01$ & $671$ & $0.05$ & $0.001$ & $1091$\\
\hline $0.01$ & $0.05$ & $3501$ & $0.01$ & $0.01$ & $6051$ & $0.01$ & $0.001$ & $9801$\\
\hline $0.005$ & $0.05$ & $7401$ & $0.005$ & $0.01$ & $12701$ & $0.005$ & $0.001$ & $20701$\\
\hline $0.001$ & $0.05$ & $38501$ & $0.001$ & $0.01$ & $66501$ & $0.001$ & $0.001$ & $108001$\\
 \hline
\end{tabular}
\end{center}
\end{table}

Finally, we would like to point out that similar characteristics of the coverage probability can be shown for
the problem of estimating a Poisson parameter or the proportion of a finite population, which allows for the
exact computation of minimum sample size. For details, see our recent papers \cite{Chen, Chen2}.

\section{Conclusion}

Determining sample size is a very important issue because samples that are too large may waste time, resources
and money, while samples that are too small may lead to inaccurate results. We have developed an exact method
for the computation of minimum sample size for the estimation of binomial parameters, which is not computational
demanding. Our sample size method permits rigorous control of statistical sampling error. Exact previously
unavailable minimum sample sizes is obtained by implementing the new method on a personal computer. Specially,
for the convenient use of practitioners, we have obtained a MATLAB data file of sample sizes for a very large
number of combinations of margin of error and confidence level, which can be available upon request. It is hoped
that our sample size method can be useful to improve the rigorousness and efficiency of statistical inference on
the very old estimation problem of binomial parameters.

\appendix

\sect{Proof of Theorem \ref{thm_abs}}

We denote the number of successes as $K = \sum_{i=1}^n X_i$.
Define
\[ C(p)  =  \Pr \li \{ \li | \f{K}{n} - p \ri | < \vep \ri \} =
\Pr \li \{ g(p) \leq K \leq h(p) \ri \}
\] where
\[
g(p) = \lf n( p - \vep) \rf + 1, \qqu h(p) = \lc n( p + \vep) \rc
- 1.
\]
It should be noted that $C(p), \; g(p)$ and $h(p)$ are actually multivariate functions of $p, \; \vep$ and $n$.
For simplicity of notations, we drop the arguments $n$ and $\vep$ throughout the proof of Theorem \ref{thm_abs}.

We need some preliminary results.

\beL \la{minus} Let $p_\ell = \f{\ell}{n} - \vep$ where $\ell \in
\bb{Z}$. Then, $h(p) = h(p_{\ell + 1}) = \ell$ for any $p \in
(p_\ell, p_{\ell +1})$. \eeL

\bpf For $p \in ( p_\ell, \; p_{\ell + 1})$, we have $0 < n \li (p
- p_\ell \ri ) < 1$ and \bee h (p)  & = & \lc n( p + \vep) \rc - 1\\
& = & \lc n( p_\ell + \vep + p - p_\ell ) \rc - 1 \\
& = &  \li \lc n \li ( \f{\ell}{n} - \vep  + \vep + p - p_\ell  \ri )  \ri \rc - 1\\
& = & \ell - 1 + \li \lc n \li (p - p_\ell \ri )  \ri \rc\\
& = & \ell\\
& = & \li \lc n \li ( \f{\ell + 1}{n} - \vep  + \vep \ri )  \ri
\rc - 1 = h(p_{\ell + 1}). \eee

\epf

\beL \la{plus} Let $p_\ell = \f{\ell}{n} + \vep$ where $\ell \in
\bb{Z}$. Then, $g(p) = g(p_{\ell})=  \ell + 1$ for any $p \in
(p_\ell, p_{\ell +1})$. \eeL

 \bpf
For $p \in \li ( p_\ell, \; p_{\ell + 1} \ri )$,  we have $-1 < n
\li (p - p_{\ell + 1} \ri ) < 0$ and \bee g (p) & = & \lf n( p -
\vep) \rf + 1\\
& = &  \lf n( p_{\ell + 1} - \vep + p - p_{\ell + 1} ) \rf + 1\\
& = &  \li \lf n \li ( \f{ \ell + 1 } { n  } + \vep - \vep \ri ) \ri
\rf +
\lf n( p - p_{\ell + 1} ) \rf +  1\\
& = & \li \lf n \li ( \f{ \ell + 1 } { n  } + \vep  - \vep \ri ) \ri \rf - 1 +  1\\
& = & \ell + 1\\
& = &  \li \lf n \li ( \f{ \ell } { n  } + \vep - \vep \ri ) \ri
\rf +  1 = g(p_\ell). \eee

\epf

\beL \la{constant} Let $\al < \ba$ be two consecutive elements of
the ascending arrangement of all distinct elements of $\{a, b \}
\cup
 \{ \f{\ell}{n} + \vep \in (a, b) : \ell \in \bb{Z} \} \cup
 \{ \f{\ell}{n} - \vep \in (a, b) : \ell \in \bb{Z} \} $.
Then, both $g(p)$ and $h (p)$ are constants for any $p \in (\al, \ba)$.
 \eeL

 \bpf
Since $\al$ and $\ba$ are two consecutive elements of the
ascending arrangement of all distinct elements of the set, it must
be true that there is no integer $\ell$ such that $\al <
\f{\ell}{n} + \vep < \ba$ or $\al < \f{\ell}{n} - \vep < \ba$.  It
follows that there exist two integers $\ell$ and $\ell^\prime$
such that {\small $(\al, \ba) \subseteq \li ( \f{\ell}{n} + \vep,
\f{\ell + 1}{n} + \vep \ri )$} and {\small $(\al, \ba) \subseteq
\li ( \f{\ell^\prime}{n} - \vep, \f{\ell^\prime + 1}{n} - \vep \ri
)$.}  Applying Lemma \ref{minus} and Lemma \ref{plus}, we have
{\small $g(p) = g \li ( \f{\ell}{n} + \vep \ri )$} and {\small
$h(p) = h \li (\f{\ell^\prime + 1}{n} - \vep \ri )$} for any $p
\in (\al, \ba)$.

 \epf

\beL \la{lem_lim}
 For any $p \in (0,1)$, $\lim_{\eta \downarrow 0} C(p + \eta) \geq C(p)$
 and $\lim_{\eta \downarrow 0} C(p - \eta) \geq C(p)$.
\eeL

\bpf Observing that $h(p + \eta) \geq h(p)$ for any $\eta > 0$ and
that
 \bee g(p + \eta) & =  &  \lf n( p + \eta - \vep)
\rf + 1 \\
& = & \lf n( p - \vep) \rf + 1 + \lf n( p - \vep) -
\lf n( p  - \vep) \rf + n \eta \rf \\
 & = & \lf n( p - \vep)
\rf + 1  = g(p) \eee for $0 < \eta < \f{ 1 + \lf n( p - \vep) \rf
-  n( p - \vep)} {n}$, we have \be \la{ineqa} S(n, g(p + \eta), h
(p + \eta), p + \eta ) \geq S(n, g(p), h (p), p + \eta ) \ee for
$0 < \eta < \f{ 1 + \lf n( p - \vep) \rf -  n( p - \vep)} {n}$.
Since
\[ h(p + \eta)  =   \lc n( p + \eta + \vep) \rc - 1 = \lc n( p +
\vep) \rc - 1 + \lc n( p + \vep) - \lc n( p + \vep) \rc + n \eta
\rc, \]  we have \[ h(p + \eta) = \bec \lc n( p + \vep) \rc &
\tx{for} \; n( p + \vep) = \lc n( p + \vep) \rc \; \tx{and} \; 0 <
\eta <
\f{1}{n},\\
\lc n( p + \vep) \rc  - 1 & \tx{for} \; n( p + \vep) \neq \lc n( p
+ \vep) \rc \; \tx{and} \; 0 < \eta < \f{\lc n( p + \vep) \rc - n(
p + \vep)}{n}. \eec
\]
It follows that both $g(p + \eta)$ and $h(p + \eta)$ are
independent of $\eta$ if $\eta > 0$ is small enough.  Since $S(n,
g, h, p + \eta)$ is continuous with respect to $\eta$ for fixed
$g$ and $h$, we have that $\lim_{\eta \downarrow 0} S(n, g(p +
\eta), h (p + \eta), p + \eta )$ exists.  As a result, \bee
\lim_{\eta \downarrow 0} C(p + \eta) & = & \lim_{\eta \downarrow
0} S(n,
g(p + \eta), h (p + \eta), p + \eta )\\
& \geq & \lim_{\eta \downarrow 0} S(n, g(p), h (p), p + \eta ) =
S(n,  g(p), h (p), p ) = C(p), \eee where the inequality follows
from (\ref{ineqa}).

Observing that $g(p - \eta) \leq g(p)$ for any $\eta > 0$ and that
 \bee h(p - \eta) & =  & \lc n( p - \eta + \vep)
\rc - 1 \\
& = & \lc n( p + \vep) \rc - 1 + \lc n( p + \vep) -
\lc n( p  + \vep) \rc - n \eta \rc\\
 & = & \lc n( p + \vep)
\rc - 1  = h(p) \eee for $0 < \eta < \f{ 1 + n( p + \vep) - \lc n(
p + \vep) \rc } {n}$, we have \be \la{ineqb} S(n, g(p - \eta), h
(p - \eta), p - \eta ) \geq S(n, g(p), h (p), p - \eta ) \ee for
$0 < \eta < \f{ 1 + n( p + \vep) - \lc n( p + \vep) \rc } {n}$.
Since \[ g(p - \eta)  =   \lf n( p - \eta - \vep) \rf + 1 = \lf n(
p - \vep) \rf + 1 + \lf n( p - \vep) - \lf n( p  - \vep) \rf - n
\eta \rf,
\]
we have
\[
g(p - \eta) = \bec \lf n( p - \vep) \rf & \tx{for} \; n( p - \vep)
= \lf n( p  - \vep) \rf \; \tx{and} \; 0 < \eta < \f{1}{n},\\
\lf n( p - \vep) \rf  + 1 & \tx{for} \; n( p - \vep) \neq \lf n( p
- \vep) \rf \; \tx{and} \; 0 < \eta < \f{n( p - \vep) - \lf n( p -
\vep) \rf }{n}. \eec
\]
It follows that both $g(p - \eta)$ and $h(p - \eta)$ are
independent of $\eta$ if $\eta > 0$ is small enough.  Since $S(n,
g, h, p - \eta)$ is continuous with respect to $\eta$ for fixed
$g$ and $h$, we have that $\lim_{\eta \downarrow 0} S(n, g(p -
\eta), h (p - \eta), p - \eta )$ exists. Hence, \bee \lim_{\eta
\downarrow 0} C(p - \eta) & =
& \lim_{\eta \downarrow 0} S(n, g(p - \eta), h (p - \eta), p -\eta )\\
& \geq & \lim_{\eta \downarrow 0} S(n, g(p), h (p), p -\eta ) =
S(n, g(p), h (p), p ) = C(p), \eee where the inequality follows
from (\ref{ineqb}).

\epf

\beL \la{unimodal} Let $0 < u < v < 1$ and $g \leq h$. Then,
\[ \min_{p \in [u, v]} S(n, g, h, p) = \min \{ S(n, g, h, u), \;
S(n, g, h, v) \}.
\]
 \eeL

 \bpf

 It can be checked that \be
\la{dif} \f{ \pa B(n, k, p) } { \pa p } = n [B(n-1, k -1, p) - B(n-1, k,p)]
 \ee for any integer $k$.  By
(\ref{dif}), it is ready to show that \be \la{ddd} \f{ \pa S(n, 0, l, p) } { \pa p } = - n B(n-1, l, p). \ee

To show the lemma, it suffices to consider $6$ cases as follows.

Case (i): $g \leq h < 0 < n$. In this case, $S(n, g, h, p) = 0$ for
any $p \in [u,v]$.

Case (ii): $0 < g \leq n \leq h $. In this case, $S(n, g, h, p)
 = S(n, g, n, p) = 1 - S(n, 0, g - 1, p)$, which is increasing in view of (\ref{ddd}).

Case (iii): $0 < n < g \leq h$. In this case, $S(n, g, h, p) = 0$
for any $p \in [u,v]$.

Case (iv): $g \leq 0 \leq  h < n$. In this case, $S(n, g, h, p)
 = S(n, 0, h, p)$, which is decreasing as a result of (\ref{ddd}).

Case (v): $g \leq 0 < n \leq h$.  In this case, $S(n, g, h, p) = 1$
for any $p \in [u,v]$.

Clearly, the lemma is true for the above five cases. \bsk

Case (vi): $0 < g \leq h < n$. By (\ref{ddd}),
 {\small \bee \f{ \pa C(p) } { \pa p }
& = & \f{ \pa S(n, 0, h, p) } { \pa p } - \f{ \pa S(n, 0, g-1, p) } { \pa p }\\
& = & \f{ n!  } { (g-1) ! (n-g)! } \; p^{g-1} (1 -p)^{n-g} - \f{
n! }
{ h ! (n-h-1)! } \; p^{h} (1 -p)^{n-h-1}\\
& = & \li [ \f{ h ! (n-h-1)!  } { (g-1) ! (n-g)! } - \li (
\f{p}{1-p} \ri )^ {h-g+1} \ri ] \f{ n! } { h ! (n-h-1)! } \;
p^{g-1} (1 -p)^{n-g} > 0 \eee} if {\small $p < 1 - \li \{ 1 + \li
[ \f{ h ! (n-h-1)! } { (g-1) ! (n-g)! } \ri ]^{\f{1}{h-g+1}} \ri
\}^{-1}$. }

From such investigation of the derivative of $C(p)$ with respective to $p$, we can see that one of the following
three cases must be true: (1) $C(\mu)$ decreases monotonically for $\mu \in [u, v]$; (2) $C(\mu)$ increases
monotonically for $\mu \in [u, v]$; (3) there exists a number $\se \in (u, v)$ such that $C(\mu)$ increases
monotonically for $\mu \in [u, \se]$ and decreases monotonically for $\mu \in (\se, v]$.  It follows that the
lemma must be true for Case (vi).

 \epf

\beL \la{inbetween}
Let $\al < \ba$ be two consecutive elements of
the ascending arrangement of all distinct elements of $\{a, b \}
\cup
 \{ \f{\ell}{n} + \vep \in (a, b) : \ell \in \bb{Z} \} \cup
 \{ \f{\ell}{n} - \vep \in (a, b) : \ell \in \bb{Z} \} $.  Then,
 $C(p) \geq \min \{ C(\al), \; C(\ba) \}$ for any $p \in (\al, \ba)$.
\eeL

\bpf

By Lemma \ref{constant}, $g(p)$ and $h(p)$ are constants for any $p
\in (\al, \ba)$. Hence, we can drop the argument and write $g(p) =
g, \; h(p) = h$ and $C(p) = S(n, g, h, p)$.

For $p \in (\al, \ba)$, define interval $[\al + \eta, \ba - \eta]$
with {\small $0 < \eta < \min \li (p - \al, \ba - p, \f{\ba -
\al}{2} \ri )$}. Then, $p \in [\al + \eta, \ba - \eta]$. By Lemma
\ref{unimodal},
\[
C(p) \geq \min_{\mu \in [\al + \eta, \ba - \eta]} C(\mu) = \min \{
C(\al + \eta), \; C(\ba - \eta) \}
\]
for {\small $0 < \eta < \min \li (p - \al, \ba - p, \f{\ba -
\al}{2} \ri )$}.  By Lemma \ref{lem_lim}, both $\lim_{\eta
\downarrow 0} C(\al + \eta)$ and $\lim_{\eta \downarrow 0}C(\ba -
\eta)$ exist and are bounded from below by $C(\al)$ and $C(\ba)$
respectively. Hence, \bee C(p) & \geq & \lim_{\eta \downarrow 0}
\; \min \{
C(\al + \eta), \; C(\ba - \eta) \}\\
& = & \min \li \{ \lim_{\eta \downarrow 0} C(\al + \eta), \;
\lim_{\eta \downarrow 0} C(\ba - \eta) \ri \} \geq \min \{ C(\al),
\; C(\ba) \} \eee for any $p \in (\al, \ba)$. \epf

\bsk

We are now in position to prove Theorem \ref{thm_abs}.   The statement about the minimum of the coverage
probability follows immediately from Lemma \ref{inbetween}.  The number of elements of the discrete set
described in Theorem \ref{thm_abs} can be calculated as follows.  Since
\[ a < \f{\ell}{n} - \vep < b \LRA 1 + \lf n(a + \vep) \rf  \leq \ell \leq \lc n(b + \vep) \rc - 1,
\]
the number of elements in $\{ \f{\ell}{n} - \vep \in (a, b) : \ell
\in \bb{Z} \}$ is \[ \lc n(b + \vep) \rc - \lf n(a + \vep) \rf - 1
 < n(b + \vep) + 1 - [ n(a + \vep) - 1  ] - 1 = n(b - a) + 1.
\]
Since \[
 a < \f{\ell}{n} + \vep  < b
 \LRA 1 + \lf n(a - \vep) \rf \leq \ell \leq \lc n(b - \vep) \rc -
 1,
 \]
the number of elements in $\{ \f{\ell}{n} + \vep \in (a, b) : \ell
\in \bb{Z} \}$ is \[ \lc n(b - \vep) \rc - \lf n(a - \vep) \rf - 1
< n(b - \vep) + 1 - [ n(a - \vep) - 1] - 1  = n(b - a) + 1.
\]
Hence, the total number of elements of $\{a, b \} \cup
 \{ \f{\ell}{n} + \vep \in (a, b) : \ell \in \bb{Z} \} \cup
 \{ \f{\ell}{n} - \vep \in (a, b) : \ell \in \bb{Z} \} $
  is less than $2n(b-a) +
4$.  This concludes the proof of Theorem \ref{thm_abs}.

 \sect{Proof of Theorem \ref{thm_com}}

Clearly, $\Pr \{ | \wh{\bs{p}}_n - p | < \vep \} = c_a$ for $p = a$ and $\Pr \{ | \wh{\bs{p}}_n - p | < \vep \}
= c_b$ for $p = b$.

For $p = \f{\ell}{n} - \vep \in \{ \f{\ell}{n} - \vep \in (a, b) : \ell \in \bb{Z} \}$, we have
\[
\lc n(p + \vep) \rc - 1 = \ell - 1, \qqu \lf n(p  - \vep) \rf + 1
= \lf \ell - 2 n \vep \rf + 1 = \ell + 1 - \lc 2 n \vep \rc,
\]
$c_-(\ell) = \Pr \{ | \wh{\bs{p}}_n - p | < \vep \}$ and $1 + \lf n(a + \vep) \rf  \leq \ell \leq \lc n(b +
\vep) \rc - 1$.

For $p = \f{\ell}{n} + \vep \in \{ \f{\ell}{n} + \vep \in (a, b) : \ell \in \bb{Z} \}$, we have
\[
\lf n(p - \vep) \rf + 1 = \ell + 1, \qqu \lc n(p  + \vep) \rc - 1
= \lc \ell + 2 n \vep \rc - 1 = \ell - 1 + \lc 2 n \vep \rc,
\]
$c_+(\ell) = \Pr \{ | \wh{\bs{p}}_n - p | < \vep \}$ and $1 + \lf n(a - \vep) \rf \leq \ell \leq \lc n(b - \vep)
\rc - 1$.

Hence, statement (I) of Theorem \ref{thm_com} can be shown by making use of the above observation and invoking
Theorem \ref{thm_abs}.

\bsk

To show statements (II) and (III), consider function $S(n,r,s,p) = \sum_{ k = r }^{ s } B(n, k, p)$ with $r \leq
s$. Applying (\ref{ddd}), we can show that {\small $\f{ \pa S(n,r,s,p) } { \pa p} = n [ B(n-1, r-1, p) - B(n-1,
s, p)]$ } and {\small \[ \f{ \pa^2 S(n,r,s,p) } { \pa p^2} = n (n-1) [ B(n-2, r-2, p)  + B(n-2, s, p) -  B(n-2,
r-1, p) - B(n-2, s-1, p)]. \]} Define $\vDe(n, \se, r,s) = S \li (n, r, s + 1, \se + \f{1}{n} \ri )  - S(n, r -
1, s, \se)$.  Then, {\small \[ \vDe(n, \se, r,s) = B \li (n, s + 1, \se + \f{1}{n} \ri ) -  B(n, r-1, \se) + S
\li (n,r,s, \se + \f{1}{n} \ri ) - S(n,r,s,\se ).
\]}
By Taylor's expansion formula, {\small \bee &   & S \li (n,r,s,
\se + \f{1}{n} \ri ) - S(n,r,s,\se ) = \li . \f{ \pa S(n,r,s,p) }
{ \pa p} \ri |_{p = \se} \times \f{1}{n} + \li. \f{ \pa^2
S(n,r,s,p) } {
\pa p^2} \ri |_{p = \ze} \times \f{1}{2 n^2}\\
&   & = B(n-1, r-1, \se) - B(n-1, s, \se)\\
&   &  \qu +  \f{ n-1 } {2 n } \li [ B(n-2,r-2, \ze) + B(n-2, s, \ze) -   B(n-2,r-1, \ze)  - B(n-2, s-1, \ze)
\ri ] \eee} where $\ze \in (\se, \se + \f{1}{n})$.  It follows that {\small \bee &
  & \vDe(n,\se, r,s)
 =   B(n-1, r-1, \se) +  B \li (n, s + 1, \se + \f{1}{n} \ri ) -  B(n, r-1, \se)  - B(n-1, s, \se)\\
&   &  \qu +  \f{ n-1 } {2 n } \li [ B(n-2,r-2, \ze) + B(n-2, s,
\ze) -   B(n-2,r-1, \ze)  - B(n-2, s-1, \ze)  \ri ]. \eee}
Differentiation of $B(n, k, p)$ with respect to $p$ shows that
$B(n, k, p)$ increases for $p \in [0, \f{k}{n}]$ and decreases for
$p \in [\f{k}{n}, 1]$.  As a result,
\[
 \min_{p \in [\se, \se + \f{1}{n}]} B(n, k,
p) = \udl{B} (n, k, \se), \qqu \max_{p \in [\se, \se + \f{1}{n}]}
p^k (1- p)^{n-k} = \ovl{B} (n, k, \se),
\]
leading to $\udl{\vDe} (n,\se, r,s) \leq \vDe(n,\se, r,s) \leq
\ovl{\vDe} (n,\se, r,s)$.

\bsk

To bound $1 - c_+(\ell-1)$ based on the bounds of $1 - c_+(\ell)$, we can use the relationship
\[
1 - c_+(\ell-1) = [1 - c_+(\ell)] + [ c_+(\ell) - c_+(\ell-1) ] =
[1 - c_+(\ell)] + \vDe(n, \se_\ell, r_\ell, s_\ell).
\]
Similarly, we can bound $1 - c_-(\ell-1)$ in terms of the bounds of $1 - c_-(\ell)$ by observing that
\[
1 - c_-(\ell-1) = [1 - c_-(\ell)] + [ c_-(\ell) - c_-(\ell-1) ] =
[1 - c_-(\ell)] + \vDe(n, \se_\ell^\prime, r_\ell^\prime,
s_\ell^\prime).
\]
This concludes the proof of Theorem \ref{thm_com}.

\sect{Proof of Theorem \ref{thm_rev} }

 Define \[ C(p)  =  \Pr \li \{ \li |
\f{K}{n} - p \ri | < \vep p \ri \} = \Pr \li \{ g(p) \leq K \leq
h(p) \ri \}
\] where
\[
g(p) = \lf n p(1 - \vep) \rf + 1, \qqu h(p) = \lc n p (1 + \vep) \rc
- 1.
\]

It should be noted that $C(p), \; g(p)$ and $h(p)$ are actually multivariate functions of $p, \; \vep$ and $n$.
For simplicity of notations, we drop the arguments $n$ and $\vep$ throughout the proof of Theorem \ref{thm_rev}.

We need some preliminary results.

\beL \la{minus_rev} Let $p_\ell = \f{\ell}{n (1 + \vep)}$ where
$\ell \in \bb{Z}$. Then, $h(p) = h(p_{\ell + 1}) = \ell$ for any
$p \in (p_\ell, p_{\ell +1})$. \eeL

\bpf For $p \in ( p_\ell, \; p_{\ell + 1})$, we have $0 < n (1 +
\vep) \li (p - p_\ell \ri ) < 1$ and \bee h (p)  & = &
\lc n p (1 + \vep) \rc - 1\\
& = & \lc n p_\ell (1 + \vep) + n(1 + \vep) (p - p_\ell) \rc - 1\\
& = & \li \lc n \li [ \f{\ell}{n}  + (1 + \vep) ( p - p_\ell)  \ri ]  \ri \rc - 1\\
& = & \ell - 1 + \li \lc n (1 +
\vep) \li (p - p_\ell \ri )  \ri \rc\\
& = & \ell\\
& = & \li \lc n \li [ \f{\ell + 1}{n (1 + \vep)}  \times (1 +
\vep) \ri ]  \ri \rc - 1 = h(p_{\ell + 1}). \eee

\epf

\beL \la{plus_rev} Let $p_\ell = \f{\ell}{n (1 - \vep)}$ where
$\ell \in \bb{Z}$. Then, $g(p) = g(p_{\ell})= \ell + 1$ for any $p
\in (p_\ell, p_{\ell +1})$. \eeL

 \bpf
For $p \in \li ( p_\ell, \; p_{\ell + 1} \ri )$,  we have $-1 < n
(1 - \vep)  \li (p - p_{\ell + 1} \ri ) < 0$ and \bee
g (p) & = &  \lf n p(1 - \vep) \rf + 1\\
& = & \lf n[ p_{\ell + 1} (1 - \vep) + (1 - \vep) (p - p_{\ell + 1}) ] \rf + 1\\
& = &  \li \lf n \times \f{ \ell + 1 } { n (1 - \vep) } \times (1 -
\vep)  \ri \rf +
\lf n (1 - \vep) ( p - p_{\ell + 1} ) \rf +  1 \\
& = & \li \lf n \times \f{ \ell + 1 } { n (1 - \vep) } \times (1 -
\vep)  \ri \rf   - 1 +  1 \\
& = & \ell + 1\\
& = & \li \lf n \times \f{ \ell } { n (1 - \vep) } \times (1 -
\vep) \ri \rf +  1 = g(p_\ell). \eee

\epf

\beL \la{constant_rev} Let $\al < \ba$ be two consecutive elements
of the ascending arrangement of all distinct elements of $\{a, b
\} \cup
 \{ \f{\ell}{n (1 - \vep)}  \in (a, b) : \ell \in \bb{Z} \} \cup
 \{ \f{\ell}{n (1 + \vep)} \in (a, b) : \ell \in \bb{Z} \} $.
Then, both $g(p)$ and $h (p)$ are constants for any $p \in (\al, \ba)$.
 \eeL

 \bpf
Since $\al$ and $\ba$ are two consecutive elements of the
ascending arrangement of all distinct elements of the set, it must
be true that there is no integer $\ell$ such that {\small $\al <
\f{\ell}{n (1 - \vep)}  < \ba$} or {\small $\al < \f{\ell}{n (1 +
\vep)} < \ba$}. It follows that there exist two integers $\ell$
and $\ell^\prime$ such that {\small $(\al, \ba) \subseteq \li (
\f{\ell}{n (1 - \vep)},  \f{\ell + 1}{n (1 - \vep)} \ri )$} and
{\small $(\al, \ba) \subseteq \li ( \f{\ell^\prime}{n (1 + \vep)},
\f{\ell^\prime + 1}{n (1 + \vep)} \ri )$}. Applying Lemma
\ref{minus_rev} and Lemma \ref{plus_rev}, we have {\small $g(p) =
g \li ( \f{\ell}{n (1 - \vep)} \ri )$} and {\small $h(p) = h \li
(\f{\ell^\prime + 1}{n (1 + \vep)} \ri )$} for any $p \in (\al,
\ba)$.

 \epf

\beL \la{lem_lim_rev}
 For any $p \in (0,1)$, $\lim_{\eta \downarrow 0} C(p + \eta) \geq C(p)$
 and $\lim_{\eta \downarrow 0} C(p - \eta) \geq C(p)$.
\eeL

\bpf

Observing that $h(p + \eta) \geq h(p)$ for any $\eta > 0$ and that
 \bee g(p + \eta) & =  &  \lf n (p + \eta) (1 - \vep)
\rf + 1  \\
& = &  \lf n p (1 - \vep ) \rf + 1 + \lf n p (1 - \vep ) -
\lf n p ( 1  - \vep ) \rf + n \eta (1 - \vep) \rf \\
 & = & \lf n p ( 1 - \vep )
\rf + 1   = g(p) \eee for $0 < \eta < \f{ 1 + \lf n p( 1 - \vep )
\rf -  n p( 1 - \vep )} {n (1 - \vep)}$, we have \be
\la{ineqa_rev} S(n, g(p + \eta), h (p + \eta), p + \eta ) \geq
S(n, g(p), h (p), p + \eta ) \ee for $0 < \eta < \f{ 1 + \lf n p(
1 - \vep ) \rf - n p( 1 - \vep )} {n (1 - \vep)}$.  Since
 \bee h(p + \eta) & =  & \lc n (p + \eta) (1 + \vep)
\rc - 1 \\
& = & \lc n p (1 + \vep) \rc - 1 + \lc n p (1 + \vep) - \lc n p (1
+ \vep) \rc + n \eta (1 + \vep) \rc, \eee we have {\small \[ h(p +
\eta) = \bec \lc n p(1 + \vep) \rc & \tx{for} \; n p( 1 + \vep) =
\lc n p( 1  + \vep) \rc \; \tx{and} \; 0 < \eta < \f{1}{n
(1 + \vep)}, \\
\lc n p( 1 + \vep) \rc  - 1  & \tx{for} \; n p( 1 + \vep) \neq \lc
n p( 1 + \vep) \rc \; \tx{and} \; 0 < \eta < \f{\lc n p( 1 + \vep)
\rc - n p( 1 + \vep)}{n (1 + \vep)}. \eec
\]}
It follows that both $g(p + \eta)$ and $h(p + \eta)$ are
independent of $\eta$ if $\eta > 0$ is small enough.  Since $S(n,
g, h, p + \eta)$ is continuous with respect to $\eta$ for fixed
$g$ and $h$, we have that $\lim_{\eta \downarrow 0} S(n, g(p +
\eta), h (p + \eta), p + \eta )$ exists.  As a result, \bee
\lim_{\eta \downarrow 0} C(p + \eta) & = & \lim_{\eta \downarrow
0} S(n,
g(p + \eta), h (p + \eta), p + \eta )\\
& \geq & \lim_{\eta \downarrow 0} S(n, g(p), h (p), p + \eta ) =
S(n,  g(p), h (p), p ) = C(p), \eee where the inequality follows
from (\ref{ineqa_rev}).

Observing that $g(p - \eta) \leq g(p)$ for any $\eta > 0$ and that
 \bee h(p - \eta) & =  &  \lc n( p - \eta) (1 + \vep)
\rc - 1 \\
& = &  \lc n p( 1 + \vep) \rc - 1 + \lc n p( 1 + \vep) -
\lc n p( 1  + \vep) \rc - n \eta ( 1 + \vep) \rc\\
 & = & \lc n p( 1 + \vep)
\rc - 1  = h(p) \eee for $0 < \eta < \f{ 1 + n p( 1 + \vep) - \lc
n p( 1 + \vep) \rc } {n ( 1 + \vep)}$, we have \be \la{ineqb_rev}
S(n, g(p - \eta), h (p - \eta), p - \eta ) \geq S(n, g(p), h (p),
p - \eta ) \ee for $0 < \eta < \f{ 1 + np( 1 + \vep) - \lc n p( 1
+ \vep) \rc } {n ( 1 + \vep)}$.  Since \bee g(p - \eta) & =  & \lf
n( p - \eta)(1 - \vep)
\rf + 1 \\
& = & \lf n p( 1 - \vep) \rf + 1 + \lf n p( 1 - \vep) - \lf n p( 1
- \vep) \rf - n \eta (1 - \vep) \rf, \eee we have {\small \[ g(p -
\eta) = \bec \lf n p( 1 - \vep) \rf & \tx{for} \; np( 1 - \vep) =
\lf np( 1  - \vep) \rf \; \tx{and} \; 0 < \eta < \f{1}{n (1 -
\vep)},\\
\lf np( 1 - \vep) \rf  + 1 & \tx{for} \; np( 1 - \vep) \neq \lf
np( 1 - \vep) \rf \; \tx{and} \; 0 < \eta < \f{np( 1 - \vep) - \lf
np( 1 - \vep) \rf }{n (1 - \vep)}. \eec
\]}
It follows that both $g(p - \eta)$ and $h(p - \eta)$ are
independent of $\eta$ if $\eta > 0$ is small enough.  Since $S(n,
g, h, p - \eta)$ is continuous with respect to $\eta$ for fixed
$g$ and $h$, we have that $\lim_{\eta \downarrow 0} S(n, g(p -
\eta), h (p - \eta), p - \eta )$ exists. Hence, \bee \lim_{\eta
\downarrow 0} C(p - \eta) & =
& \lim_{\eta \downarrow 0} S(n, g(p - \eta), h (p - \eta), p -\eta )\\
& \geq & \lim_{\eta \downarrow 0} S(n, g(p), h (p), p -\eta ) =
S(n, g(p), h (p), p ) = C(p), \eee where the inequality follows
from (\ref{ineqb_rev}).

\epf

By a similar argument as that of Lemma \ref{inbetween}, we have
\beL \la{inbetween_rev} Let $\al < \ba$ be two consecutive
elements of the ascending arrangement of all distinct elements of
$\{a, b \} \cup
 \{ \f{\ell}{n (1 - \vep)}  \in (a, b) : \ell \in \bb{Z} \} \cup
 \{ \f{\ell}{n (1 + \vep)} \in (a, b) : \ell \in \bb{Z} \} $.  Then,
 $C(p) \geq \min \{ C(\al), \; C(\ba) \}$ for any $p \in (\al, \ba)$.
\eeL

\bsk

Now we are ready to prove Theorem \ref{thm_rev}.  Clearly, the statement about the minimum of the coverage
probability follows immediately from Lemma \ref{inbetween_rev}.  It remains to calculate the number of elements
of the discrete set described in Theorem \ref{thm_rev}.  Since
\[ a < \f{\ell}{n (1 + \vep)} < b \LRA 1 + \lf n a (1 + \vep) \rf  \leq \ell \leq \lc nb (1 + \vep) \rc - 1,
\]
the number of elements in $\{ \f{\ell}{n (1 + \vep)} \in (a, b) :
\ell \in \bb{Z} \}$ is \[ \lc nb(1 + \vep) \rc - \lf na(1 + \vep)
\rf - 1  < nb(1 + \vep) + 1 - [ n a(1 + \vep) - 1  ] - 1 = n(b -
a)(1 + \vep)  + 1.
\]
Since \[
 a < \f{\ell}{n (1 - \vep)}  < b
 \LRA 1 + \lf n a(1 - \vep) \rf \leq \ell \leq \lc n b(1 - \vep) \rc -
 1,
 \]
the number of elements in $\{ \f{\ell}{n (1 - \vep)} \in (a, b) :
\ell \in \bb{Z} \}$ is \[ \lc n b (1 - \vep) \rc - \lf n a(1 -
\vep) \rf - 1 < n b(1 - \vep) + 1 - [ n a (1 - \vep) - 1] - 1  =
n(b - a) (1 - \vep) + 1.
\]
Hence, the total number of elements of $\{a, b \} \cup
 \{ \f{\ell}{n (1 - \vep)} \in (a, b) : \ell \in \bb{Z} \} \cup
 \{ \f{\ell}{n (1 + \vep)} \in (a, b) : \ell \in \bb{Z} \} $
  is less than $2n(b-a) + 4$.  The proof of Theorem \ref{thm_rev} is thus completed.

\end{document}